\documentclass[12pt]{amsart}
\usepackage{a4wide}
\usepackage{url}
\usepackage{amssymb}
\usepackage{mathbbol}

\newcommand{\charfunc}{\mathbb{1}}

\newcommand{\abs}[1]{\left\lvert #1 \right\rvert}  
\newcommand{\norm}[1]{\left\lVert #1 \right\rVert} 

\newcommand{\alg}[1]{\mathcal{#1}}
\newcommand{\hilbH}{\mathcal{H}}

\newcommand{\CC}{\mathcal{C}}
\newcommand{\cfspace}{\mathfrak{C}}

\DeclareMathOperator{\IE}{\mathbf{E}}                     

\renewcommand{\phi}{\varphi{}}

\newcommand{\IC}{\mathbf{C}}                     
\newcommand{\IN}{\mathbf{N}}                     
\newcommand{\IZ}{\mathbf{Z}}                     
\newcommand{\IP}{\mathbf{P}}                     

\makeatletter
  \DeclareRobustCommand\em
         {\@nomath\em \ifdim \fontdimen\@ne\font >\z@
                       \upshape \else \slshape \fi}
\makeatother

\newtheoremstyle{mythm}
  {9pt}
  {9pt}
  {\slshape}
  {0pt}
  {\bfseries}
  {.}
  { }
  {\thmname{#1} \thmnumber{ #2}\thmnote{ (#3)}}

\theoremstyle{mythm}

\newtheorem{Theorem}{Theorem}[section]
\newtheorem{Lemma}[Theorem]{Lemma}

\newtheorem{Corollary}[Theorem]{Corollary}

\theoremstyle{definition} 
\newtheorem{Definition}[Theorem]{Definition}
\newtheorem{Remark}[Theorem]{Remark}

\numberwithin{equation}{section}

\begin{document}
\title[Eigenspaces of Lamplighter Random Walks on Graphs]{On the Eigenspaces of
  Lamplighter Random Walks and Percolation Clusters on Graphs}
\author{Franz Lehner}
\begin{abstract}
  We show that the Plancherel measure of the lamplighter random walk on a
  graph coincides with the expected
  spectral measure of the absorbing random walk on the Bernoulli percolation
  clusters. In the subcritical regime the spectrum is pure point and we construct
  a complete orthonormal basis consisting of finitely supported eigenfunctions.
\end{abstract}


\date{\today}
\maketitle{}

\section{Introduction}
In this note we will study \emph{lamplighter random walks} on graphs.
First we show the equivalence of the lamplighter random
walk with the random walk on percolation clusters
(Theorem~\ref{thm:pnix=epCx}) by comparing moments.
The main result is a spectral resolution of its transition operator
(Theorem~\ref{thm:eigenprojections}), which shows in more detail
the intimate connection between lamplighter random walks and percolation
clusters.
The present note  complements our previous paper 
\cite{LehnerNeuhauserWoess:spectrum} 
where the case of Cayley graphs was treated using
group algebra techniques and we refer there for a more detailed discussion.
Here we point out certain simplifications which occur
when the group structure is neglected.

\section{Random Walks on Graphs}

Let $G=(V,E)$ be a graph and consider a nearest neighbour random walk
on $G$, i.e., a Markov chain $X_n$ with state space $V$ and transition
probabilities 
\begin{equation}
  \label{eq:rw:p(x,y)}
\IP[X_{n+1}=y| X_n=x] = 
  \begin{cases}
  p(x,y) &\text{if $y\sim x$} \\
  0 & \text{otherwise}
\end{cases}
\end{equation}
where $p(x,y)$ are fixed given probabilities.
Associated to this Markov chain is the \emph{transition operator}
$$
Tf(x) = \sum_{y\sim x} p(x,y)\,f(y)
$$
on the space
$$
\ell_2(G) = \bigl\{f:V\to \IC : \sum_{x\in V} \abs{f(x)}^2 < \infty\bigr\}
$$
of square summable functions on (the vertices of) $G$.
The $n$-step return probability
\begin{equation}
  \label{eq:SRWnstep}
  p^{(n)}(x,x) = \IP[X_n=x | X_0=x]
  = \sum_{x_1,x_2,\dots,x_{n-1}} p(x,x_1)\,p(x_1,x_2)\,\dotsm\, p(x_{n-1},x)
\end{equation}
can be written as
$$
p^{(n)}(x,x) = \langle T^n\delta_x,\delta_x\rangle
$$
We will assume throughout that the transition operator is selfadjoint,
and in this case the transition probabilities can be interpreted as the moments
$$
p^{(n)}(x,x) = \int t^n d\mu(t)
$$
of the \emph{Plancherel measure}
$$
d\mu(t) = \langle dE(t)\, \delta_x,\delta_x \rangle
$$
where $dE$ is comes from the spectral resolution
$$
T = \int t \,dE(t)
$$
of the operator $T$.

\section{Lamplighter random walks}

Let us  define the lamplighter random walk on the graph $G$.
We equip every vertex of $G$ with a lamp of $m$ colors or states,
one of the colors being black, i.e., the lamp turned off.
Denote $H$ the (finite) set of possible states of a lamp and $m=\abs{H}$ the 
number of different states.

\begin{Definition}
A \emph{configuration} on $G$ is a function
$\eta:G\to H$ with finite support which will be interpreted
as a state of the whole system of lamps.
We denote $\cfspace=H^{(G)}$ the set of configurations.
\end{Definition}

The \emph{switch-walk-switch lamplighter random walk} describes a random walker
moving around in the graph according to the the law \eqref{eq:rw:p(x,y)}.
Before and after each step he or she  changes the states of the lamp in the
current position at random. 
This way we obtain a Markov chain on the configuration space
$\cfspace\times G$ with transition probabilities
$$
\tilde{p}(\xi,x; \eta,y)  =  \frac{p(x,y)}{m^2}
$$
if $y\sim x$ and $\xi$ and $\eta$ coincide outside $x$ and $y$;
otherwise, $\tilde{p}(\xi,x; \eta,y)=0$.
One can interpret the configuration space $\cfspace\times G$ again as a graph
(the so-called lamplighter graph, see \cite{BartholdiWoess:2005:spectral} for a
generalization) and the $n$-step return probability 
$\tilde{p}^{(n)}(\xi,x; \eta,y)$ can be expressed by the same  
expression as \eqref{eq:SRWnstep} above.
However, because of the assumption that each color is switched with the
same probability, there is a simplification.

We start with the configuration $\iota$ where all lamps are off
and the random walker is in some start vertex $x$.
For each path $x=x_0,x_1,\dots,x_n=x$  the intermediate states
of the lamps are not important and only the last visit at each site counts,
therefore we get
\begin{multline*}
\IP[\text{at the end all lamps are off}] \\
    = \IP[\text{at the last visit each lamp is turned off}] 
    = \left(
         \frac{1}{m}
       \right)^{\abs{\{x_0,x_1,\dots,x_n\}}}
\end{multline*}
and the $n$-step return probability is
\begin{equation}
  \label{eq:LRWnstep}
  \tilde{p}^{(n)}(\iota,x;\iota,x) 
  = \sum_{x_1,x_2,\dots,x_{n-1}} p(x,x_1)\,p(x_1,x_2)\,\dotsm p(x_{n-1},x) 
    m^{-\abs{\{x_0,x_1,\dots,x_n\}}}
  .
\end{equation}

\section{Bernoulli Percolation}
Let $0<p<1$.
On the same graph $G$, consider Bernoulli site percolation with parameter $p$,
i.e., on the probability space $\Omega=\{0,1\}^G$ we consider the
independent random variables $(Y_x)_{x\in G}$ with Bernoulli distribution
$\IP[Y_x=1]=p$. 
Given $\omega\in\Omega$, let $A(\omega)$ denote the subgraph
of $G$ induced on $\{x : Y_x(\omega)=1\}$ and 
for any vertex $x\in G$, let
$C_x(\omega)$ denote the connected component of $A(\omega)$ containing a vertex $x$,
which is called the percolation cluster at $x$.
It is well known that for every graph $G$ there is a critical parameter $p_c$
such that for any vertex $x\in G$ a phase transition occurs in the sense that
for $p<p_c$ the cluster $C_x$ is almost surely finite and for
$p>p_c$ it is infinite with positive probability.
In order to make use of this fact we recall a combinatorial interpretation
of criticality.
\begin{Definition}
  For a subset $A\subseteq G$ we denote its \emph{vertex boundary}
  $$
  dA = \{y\in G : y\not\in A, y\sim x \text{ for some $x\in A$}\}
  .
  $$
  For $x\in G$, we denote
  $$
  \CC_x=\{A\subseteq G: x\in A, \text{ $A$ finite and connected}\}
  $$
  the set of finite path-connected neighbourhoods of $x$.
  In the case of $G=\IZ^2$, these are sometimes called \emph{lattice animals}.
\end{Definition}
The percolation cluster $C_x$ is finite if and only if $C_x=A$
for some $A\in\CC_x$. The latter occurs with probability
$$
\IP[C_x=A] = p^{\abs{A}} (1-p)^{\abs{dA}}
;
$$
thus for $p<p_c$ we have
\begin{equation}
\label{equ:sump1-p=1}
\sum_{A\in \CC_x} p^{\abs{A}} (1-p)^{\abs{dA}} = 1
.
\end{equation}

Now consider the absorbing random walk on $C=C_x(\omega)$:
$$
p_C(y,z) =
\begin{cases}
  p(y,z) &\text{if $y,z\in C$}\\
  0      & \text{otherwise}
\end{cases}
$$
here the $n$-step return probability is
\begin{align*}
  p_C^{(n)}(x,x) 
  &= \sum_{x_1,x_2,\dots,x_{n-1}} p_C(x,x_1)\,p_C(x_1,x_2)\,\dotsm
  p_C(x_{n-1},x) \\
  &= \sum_{x_1,x_2,\dots,x_{n-1}} p(x,x_1)\,p(x_1,x_2)\,\dotsm p(x_{n-1},x)\,
  \charfunc_{[\{x,x_1,\dots,x_{n-1}\}\subseteq C]} 
\end{align*}
Now taking the expectation of this return probability we get
$$
\IE   p_C^{(n)}(x,x)  
  = \sum_{x_1,x_2,\dots,x_{n-1}} p(x,x_1)\,p(x_1,x_2)\,\dotsm p(x_{n-1},x)\,
  p^{\abs{\{x,x_1,\dots,x_{n-1}\}}}
$$
since 
$$
\IE  \charfunc_{[\{x,x_1,\dots,x_{n-1}\}\subseteq C]} = 
\IP  [\{x,x_1,\dots,x_{n-1}\}\subseteq C] = 
  p^{\abs{\{x,x_1,\dots,x_{n-1}\}}}
.
$$
Comparing with formula~\eqref{eq:LRWnstep} we obtain the following
generalization of \cite[Theorem~1.1]{LehnerNeuhauserWoess:spectrum}.
\begin{Theorem}
\label{thm:pnix=epCx}
If we set the percolation parameter $p=1/m$, we have
$$
\tilde{p}^{(n)} (\iota,x;\iota,x) = \IE   p_C^{(n)}(x,x)
$$
and therefore the Plancherel measure of $\widetilde{T}$ coincides with the
integrated density of states of the random walk on the percolation cluster $C_x$.
\end{Theorem}

\section{Eigenfunctions}
In this section we construct eigenfunctions of the transition operator
$\widetilde{T}$ of the lamplighter random walk. 
We identify $\hilbH = l_2(\cfspace\times G)= l_2(\cfspace)\otimes_2 l_2(G)$ with the
space of square summable functions on $\cfspace\times G$.
Denote the standard bases of $l_2(H)$ and $l_2(G)$ by
$(e_\eta)_{\eta\in\cfspace}$ and $(e_x)_{x\in G}$, respectively,
and the canonical basis elements of $\hilbH$ by $e_{\eta,x} = e_{\eta}\otimes e_x$.
For each $x\in G$ we need two projections,
on the one hand the one-dimensional projections
$$
P_x e_y = \delta_{xy} e_y
$$
on $l_2(G)$ and on the other hand
the averaging operators on $l_2(\cfspace)$ given by
$$
\Theta_x f(\eta) = \frac{1}{m} \sum_{\eta' } f(\eta')
$$
where $\eta'$ runs over all $m$ different configurations which coincide with
$\eta$ outside $x$. We will denote its amplification to $\hilbH$ by
$\widetilde{\Theta}_x = \Theta_x\otimes I$.
Let us denote by $S_{xy}$ the partial isometry 
$$
S_{xy} e_{z} = \delta_{yz} e_{x}
$$
on $l_2(G)$.
Using this notation we can write the transition operator
$$
\widetilde{T}f(\eta,x) =  \sum_{y\sim x} p(x,y)\, \widetilde{\Theta}_x\widetilde{\Theta}_y f(\eta,y)
$$
as
$$
\widetilde{T}= \sum_x \sum_{y\sim x} p(x,y)\, \Theta_x \Theta_y \otimes S_{yx}
$$
where the sum converges in the strong operator topology.
Note that the operators $\widetilde{\Theta}_x$ commute with each other and with $S_{yz}$
and therefore they commute with $\widetilde{T}$.
\begin{Definition}
  For a finite subset $A\subseteq G$ with vertex boundary $dA$
  we define the projection
  $$
  \Theta_{A,dA} = \prod_{x\in A} \Theta_x \prod_{y\in dA} (I-\Theta_y)
  $$
  and its amplification to $\hilbH$ by
  $\widetilde{\Theta}_{A,dA}=\Theta_{A,dA}\otimes I$. 
\end{Definition}

The following lemma collects the main properties 
of the projections $\Theta_{A,dA}$. 
\begin{Lemma}
  \label{lem:ThetaAdA}
  \begin{enumerate}
   \item If $A$ and $B$ are connected subsets of $G$ with $A\cap
    B\ne\emptyset$,
    then $\Theta_{A,dA}\Theta_{B,dB}=0$.
   \item In the subcritical regime $p<p_c$ and for fixed $x\in G$,
    the family $(\Theta_{A,dA})_{A\in \CC_x}$ is a partition of unity on $l_2(\cfspace)$,
    i.e., different projections are mutually orthogonal and their strong sum is $I$.
  \end{enumerate}
\end{Lemma}
\begin{proof}
  If both $A\ne B$ are connected subgraphs both containing a common vertex $x$,
  then one of them must intersect the vertex boundary of the other,
  in some vertex $y$, thus $\Theta_y$ collides with its complement $I-\Theta_y$
  and they annihilate each other. Therefore $\Theta_{A,dA}$ and $\Theta_{B,dB}$
  are orthogonal.
 
  To show completeness, consider the von Neumann algebra
  $\alg{A}$ generated by $\{\Theta_x : x\in G\}$. 
  The set of vectors $\Omega=\{e_{\eta} : \eta\in \cfspace\}$ is cyclic for $\alg{A}$,
  i.e., the linear span of $\{\Theta_x e_{\eta} : x\in G, \eta\in\cfspace\}$ is dense.
  Since $\alg{A}$ is a commutative algebra, $\Omega$ is also separating.
  Thus for the proof of the lemma it suffices to show
  that $\sum_{A\in \CC_x} \Theta_{A,dA} e_{\eta} = e_{\eta}$ for each individual $\eta$.
  Since we have already shown that the projections $(\Theta_{A,dA}
  e_{\eta})_{A\in \CC_x}$ are mutually orthogonal it suffices to show that
  \begin{equation}
    \label{eq:sumnormThetaeeta}
    \sum_{A\in \CC_x} \norm{\Theta_{A,dA} e_{\eta}}^2 = 1    
  .
  \end{equation}
  Now
  \begin{equation*}
    \norm{\Theta_{A,dA} e_{\eta,x}}^2     
    = \langle \Theta_{A,dA} e_{\eta,x}, e_{\eta,x} \rangle 
    = p^{\abs{A}} (1-p)^{\abs{dA}}
    .
  \end{equation*}
  and thus condition \eqref{eq:sumnormThetaeeta} is equivalent to
  condition \eqref{equ:sump1-p=1}, which
  is satisfied in the subcritical regime (and sometimes in the critical regime
  as well).
\end{proof}
\begin{Theorem}
  \label{thm:eigenprojections}
  In the subcritical regime, 
  $( \Theta_{A,dA} \otimes P_A )_{A \subseteq G \text{ finite, connected}}$ is a partition of unity on
  $\hilbH$
  and reduces $\widetilde{T}$:
  \begin{equation}
    \label{eq:tThetaPA=ThetaTA}
    \widetilde{T} \,(\Theta_{A,dA} \otimes P_A)   = \Theta_{A,dA} \otimes T_A 
  \end{equation}
  where
  $$
  T_A = P_A T P_A
  $$
  denotes the truncation of the transition operator $T$ of the simple random  walk on $T$
  to the subgraph $A$.
\end{Theorem}
\begin{proof}
  Indeed it follows from Lemma~\ref{lem:ThetaAdA}
  that the family $(\Theta_{A,dA}\otimes{} P_x)_{x\in G, A\in \CC_x}$ is also a partition of unity:
  $$
  \sum_{x\in G} \sum_{A\in \CC_x} \Theta_{A,dA}\otimes P_x = I\otimes I
  $$
  and rearranging this sum we obtain
  $$
  \sum_A \Theta_{A,dA} \otimes P_A = I\otimes I
  .
  $$
  Now we show that each $\Theta_{A,dA}\otimes P_A$ reduces $\widetilde{T}$.
  First note that $\Theta_{A,dA}$ commutes with all $\Theta_x$ and 
  therefore $\widetilde{\Theta}_{A,dA}$ commutes with $\widetilde{T}$.
  It suffices to show equality of left and right hand side of 
  \eqref{eq:tThetaPA=ThetaTA} evaluated at $f=e_{\eta,x}$ for every
  $\eta\in{}\cfspace{}$ and $x\in{} A$.
  If $x\not\in A$ then both sides vanish and there is nothing to show. 
  Assume now that $x\in A$, then
  \begin{align*}
    \widetilde{T}\, (\Theta_A\otimes P_A)\, e_{\eta,x} &=  \widetilde{\Theta}_{A,dA} \widetilde{T} e_{\eta,x} \\
    &= \widetilde{\Theta}_{A,dA} \sum_{y\sim x} p(x,y)\, \widetilde{\Theta}_x\widetilde{\Theta}_y e_{\eta,y} \\
    &= \widetilde{\Theta}_{A,dA} \sum_{\substack{y\sim x\\ y\in A}} p(x,y)\, e_{\eta,y} \\
    &= (\Theta_{A,dA}\otimes T_A)\,e_{\eta,x}
    .
  \end{align*}
\end{proof}
\begin{Corollary}
  In the subcritical regime
  there exists a complete  orthonormal system of finitely supported
  eigenfunctions of $\widetilde{T}$. 
\end{Corollary}
\begin{proof}
  If we denote for each finite connected subgraph $A\subseteq{}G$ by 
  $\{f_a: a\in{} A\}$ a basis of $l_2(A)\subseteq{}l_2(G)$ consisting of eigenfunctions 
  of $T_A$, the eigenspaces of $\widetilde{T}$ are given by
  $\{\Theta_{A,dA} l_2(\cfspace{}) \otimes f_a : A \subseteq G\ \text{connected},
  a\in A\}$. For every finite connected subset $A\subseteq G$
  we  apply the Gram-Schmidt procedure to $\{\Theta_{A,dA} e_\eta:
  \eta\in\cfspace{}\}$ to obtain a basis of $\Theta_{A,dA}l_2(\cfspace{})$
  consisting of finitely supported functions $(\phi_{A,i})_{i\in \IN}$. This is
  possible because $\Theta_{A,dA} e_\eta$ has finite support for each
  $\eta\in\cfspace$ and $l_2(\cfspace{})$ is separable. Putting these functions
  together with the eigenfunctions $f_a$ we obtain the eigenbasis
  $$
  \{ \phi_{A,i}\otimes f_a : A\subseteq G, i\in\IN, a\in A\}
  .
  $$
\end{proof}
\begin{Remark}
  In the case where the graph $G$ is a Cayley graph as considered in
  \cite{LehnerNeuhauserWoess:spectrum}, the projections
  $\Theta_{A,dA}\otimes P_A$ are not elements of the group algebra
  and therefore provide a new partition,
  different from the one obtained in \cite{LehnerNeuhauserWoess:spectrum}.
\end{Remark}

\section{Concluding Remarks}
\begin{enumerate}
 \item 
  Similar results hold if the lamps are placed on the edges, see again
  \cite{LehnerNeuhauserWoess:spectrum}  for a discussion.
  The preceding considerations also hold 
  when one allows the number of colors (and accordingly the percolation
  parameter) to vary among the vertices $x$,
  however for the sake of simplicity  only identical lamps on
  all vertices were  considered here.
 \item 
  It is not essential that the projections $\Theta_x$ are averaging operators
  and $p=1/m$.
  As discussed in \cite{LehnerNeuhauserWoess:spectrum}, similar deterministic
  models can be constructed for arbitrary percolation parameters $p$.
 \item 
  It is still unknown what happens in the supercritical regime,
  where it is conjectured that continuous spectrum occurs at least in some cases.
  It may be hoped that the new approach will lead to some insight into this question.
\end{enumerate}

\bibliography{LamplighterGraph}
\bibliographystyle{amsplain}

\end{document}